\documentclass{amsart}
\usepackage{amscd,amssymb,verbatim}
\usepackage{amssymb,amsmath,amsthm,epsfig}

  \theoremstyle{plain}
  \newtheorem{thm}{Theorem}[section]
  \theoremstyle{definition}
  \newtheorem{defn}[thm]{Definition}
  \theoremstyle{plain}
  \newtheorem{lem}[thm]{Lemma}
  \theoremstyle{plain}
  
  \theoremstyle{plain}
  \newtheorem{prop}[thm]{Proposition}
  \theoremstyle{remark}
  \newtheorem*{claim*}{Claim}

\newcommand{\N}{\mathbb{N}}

\newcommand{\R}{\mathbb{R}}

\newcommand{\coo}{c_{00}}

\begin{document}

\title{Weakly null sequences with upper estimates}
\author{Daniel Freeman}

\address{Department of Mathematics\\ Texas A\&M University\\
College Station, TX 77843-3368} \email{freeman@math.tamu.edu}

\begin{abstract} We prove that if $(v_i)$ is a normalized basic sequence and $X$
is a Banach space such that every normalized weakly null sequence in
$X$ has a subsequence that is dominated by $(v_i)$, then there
exists a uniform constant $C\geq1$ such that every normalized weakly
null sequence in $X$ has a subsequence that is C-dominated by
$(v_i)$. This extends a result of Knaust and Odell, who proved this
for the cases in which $(v_i)$ is the standard basis for $\ell_p$ or
$c_0$.

\end{abstract}
\maketitle
\section{Introduction}\label{S1}

In some circumstances, local estimates give rise to uniform global
estimates.  An elementary example of this is that every continuous
function on a compact metric space is uniformly continuous. Uniform
estimates are especially pertinent in functional analysis, as one of
the cornerstones to the subject is the Uniform Boundedness
Principle. Because uniform estimates are always desirable, it is
important to determine when they occur.  In this paper, we are
concerned with uniform upper estimates of weakly null sequences in a
Banach space. Before stating precisely what we mean by this, we give
some historical context.

For each $1<p<\infty$, Johnson and Odell \cite{JO} have constructed
a Banach space $X$ such that every normalized weakly null sequence
in $X$ has a subsequence equivalent to the standard basis for
$\ell_p$, and yet there is no fixed $C\geq 1$ such that every
normalized weakly null sequence in $X$ has a subsequence
$C-$equivalent to the standard basis for $\ell_p$.  A basic sequence
$(x_i)$ is equivalent to the unit vector basis for $\ell_p$ if it
has both a lower and an upper $\ell_p$ estimate. That is there exist
constants $C,K\geq1$ such that:
$$\frac{1}{K}(\sum|a_i|^p)^{1/p}\leq ||\sum a_i x_i||\leq C(\sum
|a_i|^p)^{1/p}\quad \forall (a_i)\in \coo.
$$
The examples of Johnson and Odell show that the upper constant $C$
and the lower constant $K$ cannot always both be chosen uniformly.
It is somewhat surprising then that Knaust and Odell proved
\cite{KO2} that actually the upper estimate can always be chosen
uniformly.  Specifically, they proved that for every Banach space
$X$ if each normalized weakly null sequence in $X$ has a subsequence
with an upper $\ell_p$ estimate, then there exists a constant
$C\geq1$ such that each normalized weakly null sequence in $X$ has a
subsequence with a C-upper $\ell_p$ estimate. They also proved
earlier the corresponding theorem for upper $c_0$ estimates
\cite{KO1}. The standard bases for $\ell_p$, $1<p<\infty$ and $c_0$
enjoy many strong properties which Knaust and Odell employ in their
papers. It is natural to ask what are some necessary and sufficient
properties for a basic sequence to have in order to guarantee the
uniform upper estimate. In this paper we show that actually all
normalized basic sequences give uniform upper estimates.  We make
the following definition to formalize this.

\begin{defn}\label{D1.1}
Let $V=\left(v_{n}\right)_{n=1}^{\infty}$ be a normalized basic
sequence. A  Banach space X has property $\left(S_{V}\right)$ if
every  normalized weakly null sequence $\left(x_{n}\right)$ in X has
a subsequence $\left(y_{n}\right)$ such that for some constant
$C<\infty$

\begin{equation}
\label{E1.1} \qquad\qquad\left\Vert
\sum_{n=1}^{\infty}\alpha_{n}y_{n}\right\Vert \leq C\,\,\,\textrm{
for all }\left(\alpha_{n}\right)\in\coo\textrm{ with }\left\Vert
\sum_{n=1}^{\infty}\alpha_{n}v_{n}\right\Vert \leq1.
\end{equation}

X has property $\left(U_{V}\right)$ if C may be chosen uniformly. We
say that $\left(y_{n}\right)$ has a C-upper $V$-estimate (or that
$V$ C-dominates $\left(y_{n}\right)$) if (1) holds for C, and that
$\left(y_{n}\right)$ has an upper $V$-estimate (or that $V$
dominates $\left(y_{n}\right)$) if (1) holds for some C.
\end{defn}

Using these definitions, we can formulate the main theorem of our
paper as:
\begin{thm}\label{T1.2}
A Banach space has property $\left(S_{V}\right)$ if and only if it
has property $\left(U_{V}\right)$.
\end{thm}

$\left(S_{V}\right)$ and $\left(U_{V}\right)$ are isomorphic
properties of $V$, so it is sufficient to prove Theorem \ref{T1.2}
for only normalized bimonotone basic sequences.

In section 2 we present the necessary definitions and reformulate
our main results.  We break up the main proof into two parts which
we give in sections 3 and 4.  In section 5 we give some illustrative
examples which show in particular that our result is a genuine
extension of \cite{KO2} and not just a corollary.

For a Banach space $X$ we use the notation $B_X$ to mean the closed
unit ball of $X$ and $S_X$ to mean the unit sphere of $X$.  If
$F\subset X$ we denote $[F]$ to be the closed linear span of $F$ in
$X$. If $N$ is a sequence in $\N$, we denote $[N]^\omega$ to be the
set of all infinite subsequences of $N$.

This paper forms a portion of the author's doctoral dissertation, which is being prepared at Texas A\&M University
under the direction of Thomas Schlumprecht. The author thanks Dr. Schlumprecht for his invaluable help, guidance, and patience.

\section{Main Results}
Here we introduce the main definitions and theorems of the paper.
Many of our theorems and lemmas are direct generalizations of
corresponding results in \cite{KO2}.  We specify when we are able to
follow the same outline as a proof in \cite{KO2}, and also when we
are able to follow a proof exactly.

\begin{defn}
Let X be a Banach space and $V=\left(v_{n}\right)_{n=1}^{\infty}$ be
a normalized bimonotone basic sequence.  With the exception of (ii),
the following definitions are adapted from \cite{KO2}.
\begin{enumerate}
\item[(i)] A sequence $\left(x_{n}\right)$ in X is called a \textit{uV-sequence}
if $\left\Vert x_{n}\right\Vert \leq1$ for all $n\in\mathbb{N}$,
$\left(x_{n}\right)$ converges weakly to 0, and\[ \sup_{\left\Vert
\sum_{n=1}^{\infty}\alpha_{n}v_{n}\right\Vert \leq1}\left\Vert
\sum_{n=1}^{\infty}\alpha_{n}x_{n}\right\Vert <\infty.\]
$\left(x_{n}\right)$ is called a \textit{C-uV-sequence} if\[
\sup_{\left(\alpha_{n}\right)_{n=1}^{\infty}\in B_{V}}\left\Vert
\sum_{n=1}^{\infty}\alpha_{n}x_{n}\right\Vert <C.\]

\item[(ii)] A sequence $(x_n)$ in $X$ is called a \textit{hereditary
uV-sequence}, if every subsequence of $(x_n)$ is a u$V$-sequence,
and is called a \textit{hereditary C-uV-sequence} if every
subsequence of $(x_n)$ is a C-u$V$-sequence.

\item[(iii)] A sequence $\left(x_{n}\right)$ in X is called an \textit{M-bad-uV
sequence} for a constant $M<\infty$, if every subsequence
 of $\left(x_{n}\right)$ is a u$V$-sequence, and no subsequence of
$\left(x_{n}\right)$ is an M-u$V$-sequence.

\item[(iv)] An array $\left(x_{i}^{n}\right)_{i,n=1}^{\infty}$ in X is called
a \textit{bad uV-array}, if each sequence
$\left(x_{i}^{n}\right)_{i=1}^{\infty}$ is an $M_{n}$-bad
u$V$-sequence for some constants $M_{n}$ with
$M_{n}\rightarrow\infty$.
\item [(v)]$\left(y_{i}^{k}\right)_{i,k=1}^{\infty}$ is called a \textit{subarray} of
$\left(x_{i}^{n}\right)_{i,n=1}^{\infty}$, if there is a subsequence
$\left(n_{k}\right)$ of $\mathbb{N}$ such that every sequence
$\left(y_{i}^{k}\right)_{i=1}^{\infty}$ is a subsequence of
$\left(x_{i}^{n_{k}}\right)_{i=1}^{\infty}$.
\item[(vi)] A bad u$V$-array $\left(x_{i}^{n}\right)_{i,n=1}^{\infty}$ is said
to satisfy the \textit{V-array procedure}, if there exists a
subarray $\left(y_{i}^{n}\right)$ of $\left(x_{i}^{n}\right)$ and
there exists $\left(a_{n}\right)\subseteq\mathbb{R}^{+}$ with
$a_{n}\leq2^{-n}$, for all $n\in\N$, such that the weakly null
sequence $\left(y_{i}\right)$ with
$y_{i}:=\sum_{n=1}^{\infty}a_{n}y_{i}^{n}$ has no u$V$-subsequence.
\item[(vii)] X satisfies the \textit{V-array procedure} if every bad u$V$-array in X satisfies
the $V$-array procedure. X satisfies the \textit{V-array procedure
for normalized bad uV-arrays} if every normalized bad u$V$-array in
X satisfies the $V$-array procedure.
\end{enumerate}
\end{defn}
Note: A subarray of a bad u$V$-array is a bad u$V$-array. Also, a
bad u$V$-array satisfies the $V$-array procedure if and only if it
has a subarray which satisfies the $V$-array procedure.

Our Theorem \ref{T1.2} is now an easy corollary of the theorem
below.
\begin{thm}\label{T2.2}
Every Banach space satisfies the $V$-array procedure for normalized
bad u$V$-arrays.
\end{thm}
Theorem \ref{T2.2} implies Theorem \ref{T1.2} because if a Banach
space $X$ has property $S_V$ and not $U_V$ then there exists a
normalized bad u$V$-array, and the $V$-array procedure gives a
weakly null sequence in $B_X$ which is not u$V$; contradicting $X$
being $U_V$.

The proof for Theorem \ref{T2.2} will be given first for the
following special case.

\begin{prop}\label{P2.3}
Let K be a countable compact metric space. Then C(K) satisfies the
$V$-array procedure.
\end{prop}
The case of a general Banach space reduces to this special case by
the following proposition.

\begin{prop}\label{P2.4}
Let $\left(x_{i}^{n}\right)_{i,n=1}^{\infty}$ be a normalized bad
u$V$-array in a Banach space X. Then there exists a subarray
$\left(y_{i}^{n}\right)$ of $\left(x_{i}^{n}\right)$ and a countable
$w^{*}$-compact subset K of $B_{Y^*}$, where
$Y:=\left[y_{i}^{n}\right]_{i,n=1}^{\infty}$, such that
$\left(y_{i}^{n}|_{K}\right)$ is a bad u$V$-array in C(K).
\end{prop}
Theorem \ref{T2.2} is an easy consequence of Proposition \ref{P2.3}
and \ref{P2.4}. Note that Proposition \ref{P2.4} is only proved for
normalized bad u$V$-arrays. This makes the proof a little less
technical.

Before we prove anything about sub-arrays though, we need to first
consider just a single weakly null sequence. One of the many nice
properties enjoyed by the standard basis for $\ell_{p}$ which we
denote by $(e_i)$ is that $(e_i)$ is 1-spreading. This is the
property that every subsequence of $(e_{i})$ is 1-equivalent to
$(e_{i})$. Spreading is of particular importance because it implies
the following two properties which are implicitly used in
\cite{KO2}:
\begin{enumerate}
\item[(i)] If $(e_{i})$ C-dominates a sequence $(x_{i})$ then $(e_{i})$
C-dominates every subsequence of $(x_{i})$.
\item[(ii)]If a sequence $(x_{i})$ C-dominates $(e_{i})$ then $(x_{i})$ C-dominates every subsequence of
$(e_{i})$.
\end{enumerate}
Throughout the paper, we will be passing to subsequences and
subarrays, so properties (i) and (ii) would be very useful for us.
In our paper we have to get by without property (ii). On the other
hand, for a given sequence that does not have property (i), we may
use the following two results, which are both easy
 consequences of Ramsey's theorem (c.f. \cite{O}), and will be needed in
 subsequent sections.
\begin{lem}\label{L2.5}
Let $V=\left(v_{i}\right)_{i=1}^{\infty}$ be a normalized bimonotone
basic sequence. If $\left(x_{i}\right)_{i=1}^{\infty}$ is a sequence
in the unit ball of some Banach space X, such that every subsequence
of $\left(x_{i}\right)_{i=1}^{\infty}$ has a further subsequence
which is dominated by V then there exists a constant $1\leq
C<\infty$ and a subsequence $\left(y_{i}\right)_{i=1}^{\infty}$ of
$\left(x_{i}\right)_{i=1}^{\infty}$ so that every subsequence of
$\left(y_{i}\right)_{i=1}^{\infty}$ is C-dominated by V.
\end{lem}
\begin{proof}
Let $A_{n}=\left\{
\left(m_{k}\right)_{k=1}^{\infty}\in\left[\mathbb{N}\right]^{\omega}|\left(x_{m_{k}}\right)\textrm{
is }2^{n}\textrm{ dominated by V}\right\}$.

$A_{n}$ is Ramsey, thus for all $n\in\N$ there exists a sequence
$\left(m_{i}^{n}\right)_{i=1}^{\infty}=M_{n}\in\left[M_{n-1}\right]^{\omega}$
such that $\left[M_{n}\right]^{\omega}\subseteq A_{n}$ or
$\left[M_{n}\right]^{\omega}\subseteq A_{n}^{c}$. We claim that
$\left[M_{n}\right]^{\omega}\subseteq A_{n}$ for some
$n\in\mathbb{N}$, in which case we could choose
$\left(y_{i}\right)_{i=1}^{\infty}=\left(x_{m_{i}^{n}}\right)_{i=1}^{\infty}$.
Every subsequence of $\left(y_{i}\right)_{i=1}^{\infty}$ is then
$2^{n}$-dominated by V.

If our claim where false, we let
$\left(y_{n}\right)_{n=1}^{\infty}=\left(x_{m_{n}^{n}}\right)_{n=1}^{\infty}$
and $\left(y_{k_{n}}\right)_{n=1}^{\infty}$ be a subsequence of
$\left(y_{n}\right)_{n=1}^{\infty}$ for which there exists
$C<\infty$ such that $\left(y_{k_{n}}\right)_{n=1}^{\infty}$ is
C-dominated by V. Let $N\in\N$ such that $2^{N}-2N>C$ and set
$$\ell_{i}= \begin{cases}
    m_{i}^{N} &\text{if $i\leq N$},\\
    m_{k_{i}}^{k_{i}} &\text{ if $i>N$.}
          \end{cases}$$

Then
 $\left(\ell_{i}\right)_{i=1}^{\infty}\in\left[M_{N}\right]^{\omega}\subset A_{N}^{c}$  which implies that some
  $\left(a_{i}\right)_{i=1}^{L}\subset[-1,1]$ exists such that $\left\Vert
\sum_{i=1}^{L}a_{i}v_{i}\right\Vert \leq1$ and $\left\Vert
\sum_{i=1}^{L}a_{i}x_{\ell_{i}}\right\Vert >2^{N}$.  This yields

\begin{align*}
2^{N}<\left\Vert \sum_{i=1}^{L}a_{i}x_{l_{i}}\right\Vert &\leq\sum_{i=1}^{N}\left|a_{i}\right|+\left\Vert \sum_{i=N+1}^{L}a_{i}x_{m_{k_{i}}^{k_{i}}}\right\Vert \\
&\leq N+\left\Vert \sum_{i=N+1}^{L}a_{i}y_{k_{i}}\right\Vert \\
&\leq2N-\left\Vert \sum_{i=1}^{N}a_{i}y_{k_{i}}\right\Vert +\left\Vert \sum_{i=N+1}^{L}a_{i}y_{k_{i}}\right\Vert \\
&\leq2N+\left\Vert \sum_{i=1}^{L}a_{i}y_{k_{i}}\right\Vert \\
\end{align*}
which implies
\[
C<2^{N}-2N<\left\Vert \sum_{i=1}^{L}a_{i}y_{k_{i}}\right\Vert. \]
Thus $\left(y_{k_{n}}\right)_{n=1}^{\infty}$ being C-dominated by V
is contradicted.
\end{proof}

The following lemma is used for a given $(x_i)$ to find a
subsequence $(y_i)$ and a constant $C\geq1$ such that $(v_i)$
C-dominates every subsequence of $(y_i)$ and that C is approximately
minimal for every subsequence of $(y_i)$.

\begin{lem}\label{L2.6}
Let $V=\left(v_{n}\right)_{n=1}^{\infty}$ be a normalized bimonotone
basic sequence, $\left(x_{n}\right)_{n=1}^{\infty}$ be a sequence in
the unit ball of some Banach space X, and $a_{n}\nearrow\infty$ with
$a_{1}=0$. If every subsequence of
$\left(x_{n}\right)_{n=1}^{\infty}$ has a further subsequence which
is dominated by V then there exists a subsequence
$\left(y_{n}\right)_{n=1}^{\infty}$ of
$\left(x_{n}\right)_{n=1}^{\infty}$ and an $N\in\N$ such that every
subsequence of $\left(y_{n}\right)_{n=1}^{\infty}$ is
$a_{N+1}$-dominated by V but not $a_{N}$-dominated by V.
\end{lem}
\begin{proof}
By the previous lemma, we may assume by passing to a subsequence
that there exists $C<\infty$ such that every subsequence of
$\left(x_{n}\right)_{n=1}^{\infty}$ is C-dominated by V. Let
$M\in\mathbb{N}$ such that $a_{M}<C\leq a_{M+1}$. For $1\leq n\leq
M$ let
\[
A_{n}=\left\{
\left(m_{k}\right)\in\left[\mathbb{N}\right]^{\omega}|\quad
\begin{matrix}
 \left(x_{m_{k}}\right)_{k=1}^{\infty}\textrm{is }a_{n+1}\textrm{-dominated by V}\\
\textrm{ and is not }a_{n}\textrm{-dominated by V}.\\
\end{matrix}
\right\}
\]

$A_{n}$ is Ramsey, and $\left\{ A_{n}\right\}_{n=1}^{M}$ forms a
finite partition of $\left[\mathbb{N}\right]^{\omega}$ which implies
that there exists $N\leq M$ and
$\left(m_{k}\right)\in\left[\mathbb{N}\right]^{\omega}$ such that
$\left[\left(m_{k}\right)_{k=1}^{\infty}\right]^{\omega}\subset
A_{N}$. Every subsequence of
$\left(y_{n}\right):=\left(x_{m_{n}}\right)$ is $a_{N+1}$-dominated
by V and not $a_{N}$-dominated by V.
\end{proof}

\section{Proof of Proposition \ref{P2.3}}

Proposition \ref{P2.3} will be shown to follow easily from a
characterization of countable compact metric spaces along with
transfinite induction using the following result.

\begin{lem}\label{L3.1}
Let $\left(X_{n}\right)$ be a sequence of Banach spaces each
satisfying the V-array procedure. Then
$\left(\sum_{n=1}^{\infty}X_{n}\right)_{c_{0}}$ satisfies the
V-array procedure.
\end{lem}
To prove Lemma \ref{L3.1} we will need the following lemma which is
stated in \cite{KO2} for $\ell_p$ as Lemma 3.6.  The proof for
general $V$ follows the outline of its proof.

\begin{lem}\label{L3.2}
Let $\left(X_{n}\right)$ be a sequence of Banach spaces each
satisfying the V-array procedure and let $\left(x_{i}^{n}\right)$ be
a bad u$V$-array in some Banach space $X$. Suppose that for all
$m\in\mathbb{N}$ there is a bounded linear operator
$T_{m}:X\rightarrow X_{m}$ with $\left\Vert T_{m}\right\Vert \leq1$
such that $\left(T_{m}x_{i}^{m}\right)_{i=1}^{\infty}$ is an m-bad
u$V$-sequence in $X_{m}$. Then $\left(x_{i}^{n}\right)$ satisfies
the V-array procedure.
\end{lem}
\begin{proof}
We first consider Case 1: There exists $m\in\N$ and a subarray
$\left(y_{i}^{n}\right)$ of $\left(x_{i}^{n}\right)$ such that
$\left(T_{m}y_{i}^{n}\right)_{i,n=1}^{\infty}$ is a bad u$V$-array
in $X_{m}$. $\left(T_{m}y_{i}^{n}\right)_{i,n=1}^{\infty}$ satisfies
the $V$-array procedure because $X_{m}$ does. Therefore, there
exists a subarray $\left(T_{m}z_{i}^{n}\right)_{i,n=1}^{\infty}$ of
$\left(T_{m}y_{i}^{n}\right)_{i,n=1}^{\infty}$ and
$\left(a_{n}\right)\subset\R^{+}$with $a_n\leq2^{-n}$ such that
$(\sum_{n=1}^{\infty}a_{n}T_{m}z_{i}^{n})_{i=1}^\infty$ has no
u$V$-subsequence.
$\left(\sum_{n=1}^{\infty}a_{n}z_{i}^{n}\right)_{i=1}^\infty$ has no
u$V$-subsequence because $\left\Vert T_{m}\right\Vert \leq1$.
Therefore $\left(y_{i}^{n}\right)_{i,n=1}^{\infty}$ and hence
$\left(x_{i}^{n}\right)_{i,n=1}^{\infty}$ satisfies the V-array
procedure.

Case 2: If Case 1 is not satisfied then for all $m\in\N$ and every
subarray $(y_i^n)$ of $(x_i^n)$, we have that $(T_m y_i^n)$ is not a
bad u$V$-array in $X_m$. We may assume by passing to a subarray and
using Lemma \ref{L2.5} that there exists
$\left(N_{n}\right)_{n=1}^{\infty}\subset\N$ such that
\begin{equation}\label{E3.11}
 \left(x_{i}^{n}\right)_{i=1}^{\infty} \text{ is a hereditary }
N_{n}-uV-\text{sequence}\quad \text{for all }n\in\N.
\end{equation}

   By induction we choose for each $m\in\N_0$ a subarray $(z_{m,i}^n)_{i,n=1}^\infty$ of $(x^n_i)_{i,n=1}^\infty$
   and an $M_m\in\N$ so that
\begin{equation}\label{E3.1a}
(z_{m,i}^n)_{i,n=1}^\infty \text{ is a sub-array of
}(z_{m-1,i}^n)_{i,n=1}^\infty \quad\text{ if $m\geq 1$},
\end{equation}
\begin{equation}\label{E3.1b}
z_{m,i}^n = z^n_{m-1,i}\quad \text{ if $N_n\leq m$ and $i\in\N$,}
\end{equation}
\begin{equation}\label{E3.1c}
(T_m(z^n_{m,i}))_{i=1}^\infty \text{ is a hereditary
$M_m$-u$V$-sequence for all $n\in\N$ \quad if $m\geq 1$}.
\end{equation}

For $m=0$ let
$\left(z_{0,i}^{n}\right)_{i,n=1}^{\infty}=\left(x_{i}^{n}\right)_{i,n=1}^{\infty}$
. Now let $m\geq1$. For each $n\in\N$ such that $N_{n}\leq m$ let
$\left(z_{m,i}^{n}\right)_{i=1}^{\infty}=\left(z_{m-1,i}^{n}\right)_{i=1}^{\infty}$
and $K_{n}=m$. For each  $n\in\mathbb{N}$ such that $N_{n}>m$, using
Lemma \ref{L2.6}, we let $\left(z_{m,i}^{n}\right)_{i=1}^{\infty}$
be a subsequence of $\left(z_{m-1,i}^{n}\right)_{i=1}^{\infty}$ for
which there exists $K_{n}\in\mathbb{N}\cup\{0\}$ such that
$\left(T_{m}z_{m,i}^{n}\right)_{i=1}^{\infty}$ is a $K_{n}$-bad-u$V$
sequence and is also a hereditary
$\left(K_{n}+1\right)$-u$V$-sequence.
$\left(K_{n}\right)_{n=1}^{\infty}$ is bounded because otherwise we
are in Case 1. Let $M_{m}=\max_{n\in\N}K_{n}+1$.  This completes the
induction.

For all $n,i\in\N$ we have by (\ref{E3.1b}) that
$\left(z_{m,i}^{n}\right)_{m=1}^{\infty}$ is eventually constant.
Let
$\left(z_{i}^{n}\right)_{i,n=1}^{\infty}=\lim_{m\rightarrow\infty}\left(z_{m,i}^{n}\right)_{i,n=1}^{\infty}$.
By (\ref{E3.1c}), $\left(z_{i}^{n}\right)_{i,n=1}^{\infty}$
satisfies
\begin{equation}\label{E3.1d}
(T_m(z^n_{i}))_{i=1}^\infty \text{ is a hereditary
$M_m$-u$V$-sequence for all $m,n\in\N$}.
\end{equation}

We will now inductively choose
$\left(m_{n}\right)\in\left[\N\right]^{\omega}$ and
$\left(a_{n}\right)\subset\R^{+}$ so that for all $n\in\N$ we have:
\begin{equation}\label{E3.10}
\left(T_{m_{n}}z_{i}^{m_{n}}\right)_{i=1}^{\infty}\textrm{ is an }
m_{n}\textrm{-bad uV sequence in }X_{m_{n}}, \end{equation}
\begin{equation}\label{E3.12}a_{n}m_{n}>n,
\end{equation}
\begin{equation}\label{E3.13}\sum_{j=1}^{n-1}a_{j}N_{m_j}<\frac{a_{n}m_{n}}{4},\textrm{
and}
\end{equation}
\begin{equation}\label{E3.14.1}0<a_n<\min_{1\leq
k<n}\{2^{-n},2^{-n}\frac{a_k m_k}{4M_{m_k}}\}.
\end{equation}

Property (\ref{E3.10}) has been assumed in the statement of the
Lemma.  For n=1 let $a_{1}=\frac{1}{2}$ and $m_{1}\in\mathbb{N}$
such that $a_{1}m_{1}>1$, so (\ref{E3.12}) is satisfied.
(\ref{E3.13}) and (\ref{E3.14.1}) are vacuously true for n=1, so all
conditions are satisfied for $n=1$.

Let $n>1$ and assume $\left(a_{j}\right)_{j=1}^{n-1}$ and
$\left(m_{j}\right)_{j=1}^{n-1}$ have been chosen to satisfy
(\ref{E3.12}), (\ref{E3.13}) and (\ref{E3.14.1}). Choose $a_{n}>0$
small enough such that $a_{n}<\min_{1\leq k<n}\left\{
2^{-n},2^{-n}\frac{a_{k}m_{k}}{4M_{k}}\right\} $, thus satisfying
(\ref{E3.14.1}). Choose $m_{n}>0$ large enough to satisfy
(\ref{E3.12}) and (\ref{E3.13}). This completes the induction.

By (\ref{E3.14.1}), we have for all $n\in\N$ that
\begin{equation}\label{E3.14}
\sum_{j=n+1}^{\infty}a_{j}M_{m_{n}}<\frac{a_{n}m_{n}}{4}.
\end{equation}

We have by (\ref{E3.14.1}) that $a_j<2^{-j}$ for all $j\in\N$, so
$y_{k}:=\sum_{j=1}^{\infty}a_{j}x_{k}^{m_{j}}$ is a valid choice for
the V-array procedure. Let $C>0$ and $\left(y_{k_{i}}\right)$ be a
subsequence of $\left(y_{k}\right)$. We need to show that
$\left(y_{k_{i}}\right)$ is not a C-u$V$-sequence. Using
(\ref{E3.12}), choose $n\in\N$ so that $a_{n}m_{n}>2C$. Using
(\ref{E3.10}) choose $\ell\in\N$ and
$\left(\beta_{i}\right)_{i=1}^{\ell}\in
B_{\left[v_{i}\right]_{i=1}^{\ell}}$ such that
\begin{equation}
\label{E3.16}\left\Vert
\sum_{i=1}^{\ell}\beta_{i}T_{m_{n}}\left(x_{k_{i}}^{m_{n}}\right)\right\Vert
>m_{n}.
\end{equation}
We now have the following
\begin{align*}
\left\Vert \sum_{i=1}^{\ell}\beta_{i}y_{k_{i}}\right\Vert
&=\left\Vert
\sum_{i=1}^{\ell}\sum_{j=1}^{\infty}\beta_{i}a_{j}x_{k_{i}}^{m_j}\right\Vert \\
&\geq\left\Vert
\sum_{i=1}^{\ell}\sum_{j=n}^{\infty}T_{m_{n}}\left(\beta_{i}a_{j}x_{k_{i}}^{m_j}\right)\right\Vert
-\left\Vert
\sum_{i=1}^{\ell}\sum_{j=1}^{n-1}\beta_{i}a_{j}x_{k_{i}}^{m_j}\right\Vert
\quad\textrm{since }\left\Vert T_{m_{n}}\right\Vert \leq1 \\
&\geq
a_{n}\left\Vert
\sum_{i=1}^{\ell}\beta_{i}T_{m_{n}}x_{k_{i}}^{m_n}\right\Vert
-\sum_{j=n+1}^{\infty}a_{j}\left\Vert
\sum_{i=1}^{\ell}\beta_{i}T_{m_{n}}x_{k_{i}}^{m_j}\right\Vert
-\sum_{j=1}^{n-1}a_{j}\left\Vert
\sum_{i=1}^{\ell}\beta_{i}x_{k_{i}}^{m_j}\right\Vert\\
&>
a_{n}m_{n}-\sum_{j=n+1}^{\infty}a_{j}M_{m_{n}}-\sum_{j=1}^{n-1}a_{j}N_{m_j}\qquad\textrm{by
} (\ref{E3.16}), (\ref{E3.1d}),\textrm{and }(\ref{E3.11})\\
&\geq a_{n}m_{n}-a_{n}m_{n}/4-a_{n}m_{n}/4\qquad\textrm{by
(\ref{E3.13}) and (\ref{E3.14})} \\
&=a_{n}m_{n}/2>C. \\
\end{align*}

Therefore, $\left(y_{k_{i}}\right)$ is not a C-u$V$-sequence.
$\left(y_{i}\right)_{i=1}^{\infty}=$$\left(\sum_{j=1}^{\infty}a_{j}x_{i}^{m_{j}}\right)_{i=1}^{\infty}$
has no u$V$-subsequence, so $\left(x_{i}^{n}\right)$ satisfies the
V-array procedure which proves the lemma.
\end{proof}
Now we are prepared to give a proof of Lemma \ref{L3.1}. We follow
the outline of the proof of Lemma 3.5 in \cite{KO2}.

\begin{proof}[Proof of Lemma \ref{L3.1}]
let $\left(x_{i}^{n}\right)$ be a bad u$V$-array in $X=\left(\sum
X_{n}\right)_{c_{0}}$ and $R_{m}:X\rightarrow X_{m}$ be the natural
projections.

Claim: For all $M<\infty$ there exists $n,m\in\N$ and subsequence
$\left(y_{i}\right)_{i=1}^{\infty}$ of
$\left(x_{i}^{n}\right)_{i=1}^{\infty}$ such that
$\left(R_{m}y_{i}\right)_{i=1}^{\infty}$ is an M-bad u$V$-sequence.

Assuming the claim, we can find
$\left(N_{n}\right)_{n=1}^{\infty}\in\left[\mathbb{N}\right]^{\omega}$,
$\left(m(n)\right)_{n=1}^{\infty}\subset\N$, and subsequences
$\left(y_{i}^{n}\right)_{i=1}^{\infty}$ of
$\left(x_{i}^{N_{n}}\right)_{i=1}^{\infty}$ such that
$\left(R_{m(n)}y_{i}^{n}\right)_{i=1}^{\infty}$ is an n-bad u$V$
sequence for all $n\in\mathbb{N}$. By passing to a subsequence, we
may assume either that $m(n)=m$ is constant, or that
$\left(m(n)\right)_{n=1}^{\infty}\in\left[\mathbb{N}\right]^{\omega}$.
 If $m(n)=m$, then $R_m(y^n_i)_{n,i=1}^\infty$ is a bad u$V$-array in $X_m$.
 $R_m(y^n_i)_{n,i=1}^\infty$ satisfies the V-array procedure, and thus
 $(y^n_i)_{n,i=1}^\infty$ satisfies the V-array procedure. If
$\left(m(n)\right)_{n=1}^{\infty}\in\left[\mathbb{N}\right]^{\omega}$
let $T_{n}:=R_{m(n)}|_{\left[y_{i}^{r}\right]_{i,r=1}^{\infty}}$ and
apply Lemma \ref{L3.2} to the array
$\left(y_{i}^{n}\right)_{i,n=1}^{\infty}$ to finish the proof.

To prove the claim, we assume it is false. There exists $M<\infty$
such that for all $m,n\in\N$ every subsequence of
$\left(x_{i}^{n}\right)_{i=1}^{\infty}$ contains a further
subsequence $\left(y_{i}\right)_{i=1}^{\infty}$ such that
$\left(R_{m}y_{i}\right)_{i=1}^{\infty}$ is an M-u$V$-sequence.

By Ramsey's theorem, for each $n\in\N$ and $m\in\N$ every
subsequence of $\left(x_{i}^{n}\right)_{i=1}^{\infty}$ contains a
further subsequence $\left(y_{i}\right)_{i=1}^{\infty}$ such that
$\left(R_{m}y_{i}\right)_{i=1}^{\infty}$ is a hereditary
M-u$V$-sequence. Fix $n\in\mathbb{N}$ such that
$\left(x_{i}^{n}\right)_{i=1}^{\infty}$ is an (M+3)-bad
u$V$-sequence. We now construct a nested collection of subsequences
$\{\left(y_{k,i}\right)_{i=1}^{\infty}\}_{k=0}^\infty$ of
$\left(x_{i}^{n}\right)_{i=1}^{\infty}$ (where
$(y_{0,i})_{i=1}^\infty=(x_i^n)_{i=1}^\infty$) as well as
$\left(m_{i}\right)\in\left[\mathbb{N}\right]^{\omega}$ so that for
all $k\in\mathbb{N}$ we have

\begin{equation}\label{E3.17}
\sup_{m>m_{k}}\left\Vert R_{m}y_{k-1,k}\right\Vert \leq2^{-k},
\end{equation}

\begin{equation}\label{E3.17.2}
(y_{k,i})_{i=1}^\infty\textrm{ is a subsequence of
}(y_{k-1,i})_{i=2}^\infty,
\end{equation}

\begin{equation}\label{E3.18.1}
(R_m y_{k,i})_{i=1}^\infty \textrm{ is a hereditary M-uV-sequence }
\forall m\leq m_k.
\end{equation}

For k=1 we choose $m_{1}\in\mathbb{N}$ such that
$\sup_{m>m_{1}}\left\Vert R_{m}y_{0,1}\right\Vert \leq2^{-1}$. Pass
to a subsequence $\left(y_{1,i}\right)_{i=1}^{\infty}$ of
$\left(y_{0,i}\right)_{i=2}^{\infty}$ such that
$\left(R_{m}y_{1,i}\right)_{i=1}^{\infty}$ is a hereditary
M-u$V$-sequence for all $m\leq m_{1}$.

For $k>1$ given $m_{k-1}\in\N$ and a sequence
$\left(y_{k-1,i}\right)_{i=k}^{\infty}$. Choose $m_{k}>m_{k-1}$ so
that $\sup_{m>m_{k}}\left\Vert R_{m}y_{k-1,k}\right\Vert
\leq2^{-k}$, thus satisfying (\ref{E3.17}). Let
$\left(y_{k,i}\right)_{i=1}^{\infty}$ be a subsequence of
$\left(y_{k-1,i}\right)_{i=2}^{\infty}$ so that
$\left(R_{m}y_{k,i}\right)_{i=1}^{\infty}$ is a hereditary
M-u$V$-sequence for all $m\leq m_{k}$, thus satisfying
(\ref{E3.17.2}) and (\ref{E3.18.1}). This completes the induction.

We define $y_k=y_{k-1,k}$ for all $k\in\N$. By (\ref{E3.17.2}), we have that
$\left(y_{k,i}\right)_{i=1}^{k}\cup\left(y_{i}\right)_{i=k+1}^{\infty}$
is a subsequence of $\left(y_{k,i}\right)_{i=1}^{\infty}$.
Therefore, (\ref{E3.18.1}) gives that
\begin{equation}\label{E3.18}
(v_i)_{i=k+1}^\infty \textrm{ M-dominates
}(R_{m}y_{q_i})_{i=k+1}^{\infty}\forall m\leq
m_{k},\;\left(q_{i}\right)\in\left[\mathbb{N}\right]^{\omega},\textrm{
and }k\in\N.
\end{equation}

$\left(x_{i}^{n}\right)_{i=1}^{\infty}$ is a $(M+3)$-bad u$V$
sequence, so there exists $\left(\alpha_{i}\right)\in
B_{\left[V\right]}$ such that

\begin{equation}\label{E3.50}
\left\Vert \sum_{i=1}^{\infty}\alpha_{i}y_{i}\right\Vert
>M+3.
\end{equation}

For all $k\in\mathbb{N}$ and $m\in\left(m_{i-1},m_{i}\right]$ (with
$m{}_{0}=0$) we have that

\begin{align*}
\left\Vert
\sum_{i=1}^{\infty}R_{m}\left(\alpha_{i}y_{i}\right)\right\Vert
&\leq\sum_{i=1}^{k-1}\left|\alpha_{i}\right|\left\Vert
R_{m}y_{i}\right\Vert +\left\Vert
R_{m}\left(\alpha_{k}y_{k}\right)\right\Vert +\left\Vert
\sum_{i=k+1}^{\infty}R_{m}\left(\alpha_{i}y_{i}\right)\right\Vert \\
&\leq\sum_{i=1}^{k-1}2^{-i}+1+\left\Vert
\sum_{i=k+1}^{\infty}\alpha_{i}R_{m}(y_{i})\right\Vert
\qquad\textrm{by }(\ref{E3.17}) \\
&\leq1+1+M\qquad\textrm{by }(\ref{E3.18}) \\
\end{align*}
which implies
\[\left\Vert \sum_{i=1}^{\infty}\alpha_{i}y_{i}\right\Vert=
\sup_{m\in\N}\Vert
\sum_{i=1}^{\infty}R_{m}\left(\alpha_{i}y_{i}\right)\Vert\leq M+2.\]

 This contradicts (\ref{E3.50}), so the claim and hence the lemma
is proved.
\end{proof}
The proof for proposition \ref{P2.3} now follows in exactly the same
way as in \cite{KO2}.

\begin{proof}[Proof of Proposition \ref{P2.3}]
For every countable limit ordinal $\alpha$ we can find a sequence of
ordinals $\beta_{n}<\alpha,\:\beta_{n}\nearrow\alpha$ such that
$C(\alpha)$ is isomorphic to $\left(\sum
C\left(\beta_{n}\right)\right)_{c_{0}}$. Using induction and Lemma
\ref{L3.1} we obtain that all $C(\alpha)$-spaces, where $\alpha$ is
a countable limit ordinal, satisfy the V-array procedure. Thus, in
view of the isomorphic classification of $C(K)$-spaces for countable
compact metric spaces K (see \cite{BP}), all $C(K)$-spaces for
countable compact metric spaces K satisfy the V-array procedure.
\end{proof}

\section{Proof of Proposition \ref{P2.4}}
The proof of Theorem \ref{T2.2} will be complete once we have proven
proposition \ref{P2.4}. To make notation easier, we now consider the
triangulated version $\left(x_{i}^{n}\right)_{1\leq n\leq i<\infty}$
of the square array $\left(x_{i}^{n}\right)_{i,n=1}^{\infty}$.

\begin{lem}\label{L4.1}
A square array satisfies the $V$-array procedure if and only if it's
triangulated version does.
\end{lem}
\begin{proof}If
$\left(y_{i}^{n}\right)_{i,n=1}^{\infty}$ is a subarray of
$\left(x_{i}^{n}\right)_{i,n=1}^{\infty}$  then
$\left(y_{i}^{n}\right)_{1\leq n\leq i<\infty}$ is a triangular
subarray of $\left(x_{i}^{n}\right)_{1\leq n\leq i<\infty}$. Also if
$\left(y_{i}^{n}\right)_{1\leq n\leq i<\infty}$ is a triangular
subarray of $\left(x_{i}^{n}\right)_{1\leq n\leq i<\infty}$ then
$\left(y_{i}^{n}\right)_{1\leq n\leq i<\infty}$ may be extended to a
subarray of $\left(x_{i}^{n}\right)_{i,n=1}^{\infty}$ by letting
$\left(y_i^n\right)_{i<n}=\left(x_i^n\right)_{i<n}$.

We now show that applying the $V$-array procedure to
$\left(y_{i}^{n}\right)_{i,n=1}^{\infty}$ and
$\left(y_{i}^{n}\right)_{1\leq n\leq i<\infty}$ yield equivalent
sequences. For all $n\in\N$ let
$0\leq\left|\alpha_{n}\right|\leq2^{-n}$,
$z_{i}=\sum_{n=1}^{i}\alpha_{n}y_{i}^{n}$, and
$y_{i}=\sum_{n=1}^{\infty}\alpha_{n}y_{i}^{n}$.  For all $m\in\N$ if
$\left(\beta_{i}\right)_{i=1}^{\infty}\in B_{\left[V\right]}$ then

\begin{align*}
\left\Vert
\sum_{i=1}^{m}\beta_{i}z_{i}-\sum_{i=1}^{m}\beta_{i}y_{i}\right\Vert
&=\left\Vert
\sum_{i=1}^{m}\beta_{i}\sum_{n=i+1}^{\infty}\alpha_{n}y_{i}^{n}\right\Vert
\leq\sum_{i=1}^{m}\left|\beta_{i}\right|\sum_{n=i+1}^{\infty}\left|\alpha_{n}\right|
\leq\sum_{i=1}^{m}2^{-i}<1. \\
\end{align*}
Thus we have that $\lim_{m\rightarrow\infty}\left\Vert
\sum_{i=1}^{m}\beta_{i}z_{i}\right\Vert =\infty$ if and only if
$\lim_{m\rightarrow\infty}\left\Vert
\sum_{i=1}^{m}\beta_{i}y_{i}\right\Vert =\infty$, which implies the
claim.

\end{proof}
\begin{lem}\label{L4.2}
For all $\epsilon>0$, a triangular bad u$V$-array
$\left(x_{i}^{n}\right)_{n\leq i}$ admits a triangular subarray
$\left(y_{i}^{n}\right)_{n\leq i}$ which is basic in its
lexicographical order(where i is the first letter and n is the
second letter), and its basis constant is not greater than
$2+\epsilon$ (meaning the supremum of the norm of the projections
onto the span generated by intervals of the basis).  In other words
$y_{1}^{1},y_{2}^{1},y_{2}^{2},y_{3}^{1},y_{3}^{2},y_{3}^{3},y_{4}^{1}...$
is a basic sequence.
\end{lem}
\begin{proof}
The proof is the same as the proof that a weakly null sequence has a
basic subsequence.
\end{proof}

We now assume that the given bad u$V$-array $\left(x_{i}^{n}\right)$
is labeled triangularly and that it is a bimonotone basic sequence
in its lexicographical order. This assumption is valid because the
properties "being a bad u$V$-array" and "satisfying the V-array
procedure" are invariant under isomorphisms. We also assume that
$\left(x_{i}^{n}\right)$ is normalized.

The following theorem is our main tool used to construct the
subarray $(y^n_i)$ of $(x^n_i)$ and the countable $w^*$-compact set
$K\subset B_{[y^n_i]}$ for Proposition \ref{P2.4}.

\begin{thm}\label{T4.3} Assume that $(x^n_i)_{1\leq n\leq i}$ is a normalized
triangular array in $X$, such that for every $n\in\N$ the sequence
$(x^n_i)_{i=1}^\infty$ is weakly converging to $0$. Let $V=(v_i)$
 be a normalized basic sequence and let $(C_n)\subset [0,\infty)$
 and $\epsilon>0$.

Then $(x^n_i)$ has a triangular sub-array $(y^n_i)$ with the
following property:

For all $n,q\in\N$ and all $n\leq i_1<i_2\ldots<i_q$ all
$(\alpha_j)_{j=1}^q\in B_V$ with $\|\sum_{j=1}^q \alpha_j
y^n_{i_j}\|\geq C_n$ there is a $g\in (2+\epsilon)B_{X^*}$ and
$(\beta_j)_{j=1}^q\in B_V$,
 so that

\begin{align}
&\sum_{j=1}^q g(\beta_j y^n_{i_j})\geq C_n\label{E1},\\
&g(y^m_i)=0\text{ whenever $m\leq i$ and $i\not\in\{i_1,i_2,\ldots
i_q\}$}.\label{E2}
\end{align}

If we also assume that $(x_j^n)_{1\leq n\leq j}$ is a bimonotone
basic sequence in its lexicographical order then there exists
$(i_j)\in[\N]^\omega$ so that we may choose the sub-array $(y_j^n)$
by setting $y_j^n=x_{i_j}^n$ for all $n\leq j$.  In this case we
have the above conclusion for some $g\in(1+\epsilon)B_{Y^*}$.
\end{thm}

\begin{proof}
 After passing to a sub-array using Lemma \ref{L4.2} we can assume that $(x_i^n)$ is
 a basic sequence in its lexicographical order and that it's basis constant
 does not exceed the value $2+\epsilon/4$.
 We first renorm  $Z=[x^n_i]$  by a norm $|||\cdot|||$
 in the standard
 way so that $\|z\|\leq |||z|||\leq (2+\epsilon/4) \|z\|$ and so that
 $(x^n_i)$ is bimonotone in $Z$. We therefore can assume that $(x^n_i)$ is a bimonotone
 basis and need to show the claim of Theorem \ref{T4.3} for
 $(1+\epsilon)B_{X^*}$ instead of $(2+\epsilon)B_{X^*}$.

 Let $(\epsilon_k)\subset(0,1)$ with $\sum_{k=1}^\infty k\epsilon_k<\epsilon/4$.
 By induction on $k\in\N_0$ we choose $i_k\in\N$ and a sequence $L_k\in[\N]^\omega$,
  and  define $y^m_j=x^m_{i_j}$ for $m\leq k$ and $m\leq j\leq k$ so
  that the following conditions are satisfied.

\begin{enumerate}
\item[a)] $i_k=\min L_{k-1}<\min L_k$ and $L_k\subset  L_{k-1}$, if $k\geq 1$
($L_0=\N$).
\item[b)] For all $s,t\in\N_0$, all $1\leq m \leq k$, all
 $m\leq m_1<m_2<\ldots m_s\leq k$ and $\ell_0<\ell_1<\ldots \ell_t$ in
 $L_k$, if
 \begin{align}\label{E:1.1}
 &\exists f\in B_{X^*}
      \textrm{ with }\sum_{j=1}^{s} \alpha_j f(y_{m_j}^m)+\sum_{j=1}^t \alpha_{j+s}
 f(x^m_{\ell_j})\geq C_m\textrm{ for some } (\alpha_j)_{j=1}^{s+t}\in B_{[V]}\\
\intertext{then} \label{E:1.2}
 &\exists g\in B_{X^*}\textrm{ such that }\\
&\qquad(a)\,\sum_{j=1}^{s} \beta_j g(y_{m_j}^m)+\sum_{j=1}^t
   \beta_{j+s}g( x_{\ell_j}^m)\geq C_m
   \textrm{ for some } (\beta_j)_{j=1}^{s+t}\in B_{[V]}\notag,\\
 &\qquad(b)\,|g(y_j^{m'})|<\epsilon_j\text{ if $m'\leq k$ and $j\!\in\!\{m',\ldots
 k\}\setminus\{m_1,\ldots m_s\}$}\notag\textrm{, and}\\
 &\qquad(c)\,|g(x_{\ell_0}^{m'})|<\epsilon_{k+1} \text{ if $m'\leq k+1$.} \notag
 \end{align}
 (in the case that $s=0$ condition (b) is defined to be vacuous, also note
  that in (c) we allow $m'=k+1$).
 \end{enumerate}

We first note for $(i_j)\in[\N]^\omega$ that $(x^n_{i_j})_{n\leq j}$
is a subsequence of $(x^n_j)_{n\leq j}$ in their lexicographic
orders.  Thus $(x^n_{i_j})_{n\leq j}$ is a bimonotone basic sequence
in its lexicographic order.

For $k=0$, if $f\in B_{X^*}$ satisfies (\ref{E:1.1}) then
$g=P^*_{[x_{\ell_1}^n,\infty)}f$ satisfies (\ref{E:1.2}) by our
assumed assumed bimonotonicity.

Assume $k\geq 1$ and we have chosen $i_1<i_2<\ldots<i_{k-1}$. We let
$i_k=\min L_{k-1}$.

Fix an infinite  $M\subset L_{k-1}\setminus \{i_k\}$, a positive
integer $m\leq k$,
 an integer $0\leq s\leq k-m+1$,
and
 positive integers $m\leq m_1<m_2<\ldots m_s\leq k$ and
 define
 $$A=A(m,s,(m_j)_{j=1}^s)=\bigcap_{t\in\N_0} A_t, \text{ where }$$
$$A_t=\left\{ (\ell_j)_{j=0}^\infty\in [M]^\omega:
\begin{matrix}
 \text{ If $(m_j)_{j=1}^s$ and $(\ell_j)_{j=0}^t$ satisfy
 \eqref{E:1.1}}\\
\text{then they also satisfy \eqref{E:1.2}}
\end{matrix}
\right\}.$$

For $t\in\N$ the set $A_t$ is closed as a subset of $2^\N$ in the
product topology, thus $A$ is closed and, thus, Ramsey. We will show
that there is an infinite $L\subset M$ so that $[L]^\omega\subset
A$. Once we verified  that claim we can finish our induction step by
applying that argument successively to all choices
 of $m\leq k$,
 an $0\leq s\leq k$
and $m\leq m_1<m_2<\ldots m_s\leq k$, as there are only finitely
many.

Assume our claim is wrong and, using Ramsey's Theorem, we could find
 an $L=(\ell_j)_{j=1}^\infty$ so that
  $[L]^\omega\cap A=\emptyset$.

Let $n\in\N$ be fixed, and let $p\in\{1,2\ldots, n\} $. Then
$L^{(p)} =\{\ell_p,\ell_{n+1},\ldots \}$ is not in $A$ and
 we can choose $t_n\in\N_0$, $(\alpha_j^n)_{j=1}^{t_n+s}$ and
  $f_n\in B_{X^*}$ so that \eqref{E:1.1} is satisfied
   (for $(\ell_{n+1},\ldots \ell_{\ell+t})$
 replacing $(\ell_{1},\ldots \ell_{t})$
  and $\ell_p$ replacing $\ell_0$) but for no
    $g\in B_{X^*}$ and $(\beta_j)_{j=1}^{s+t_n}\in B_{[V]}$
    condition \eqref{E:1.2} holds.  By choosing $t_n$ to be minimal so that $\eqref{E:1.1}$ is satisfied,
    we can have $t_n$, $(\alpha_j^n)_{j=1}^{t_n+s}$ and
  $f_n$ be independent of $p$.

We now show that there is a $g_n\in B_X$ satisfying (a) and (b) of
\eqref{E:1.2}.

Let $k'=\max\big\{m-1\leq i\leq k: i\not\in\{m_1,m_2,\ldots
m_s\}\big\}$.  If $k'\leq m$ then
$\{m_1,\ldots,m_s\}=\{k'+1,k'+2,\ldots,k\}$ and by our assumed
bimonotonicity $g_n:=P^*_{[y_{k'+1}^{m},\infty)}f_n\in{B_X^*}$
satisfies (a) and (b) of \eqref{E:1.1}.  If $k'> m$ let $0\leq
s'\leq s$, such that
 $m_1<m_2<\ldots<m_{s'}<k'$, and apply the $k'-1$ step of the induction hypothesis
  to $f_n$, $(\alpha_j^n)_{j=1}^{t_n+s}$, $m\leq m_1<\ldots<m_{s'}$ (replacing $m\leq m_1<\ldots<m_s$),
 and $k'<k'+1<\ldots<m_{s}<\ell_{n+1}<\ldots<\ell_{t_n}$ (replacing $\ell_p<\ell_{n+1}<\ldots<\ell_{t_n}$)
  to obtain  a functional $g_n\in B_{X^*}$ which satisfies (a) and (b)  of \eqref{E:1.2}.

 Since $g_n$ cannot satisfy all three conditions of \eqref{E:1.2} (for
  any choice of $1\leq p\leq n$), we
  deduce that $|g_n(x_{\ell_p}^{m_p})|\geq\epsilon_{k+1}$ for some choice of
  $m_p\in\{1,2,\ldots k+1\}$.

Let $g$ be a $w^*$ cluster point of $(g_n)_{n\in\N}$. As the set
$\{1,2,\dots k+1\}$ is finite, we have for all $p\in\N_0$ that
$|g(x_{\ell_p}^{m_p})|\geq\epsilon_{k+1}$ for some
$m_p\in\{1,2,\dots k+1\}$. Which implies there exists $1\leq m\leq
k+1$ such that $|g(x^{m}_{\ell_p})|\geq\epsilon$ for infinitely many
$p\in\N$. This is a contradiction with the sequence
$(x^m_{\ell_i})_{i=1}^\infty$ being weakly null. Our claim is
verified, and we are able to fulfill the induction hypothesis.

The conclusion of our theorem now follows by the following
perturbation argument. If we have $n\leq i_1<i_2\ldots<i_q$ and
$(\alpha_j)_{j=1}^q\in B_V$ with $\|\sum_{j=1}^q \alpha_j
y^n_{i_j}\|\geq C_n$, then there exists $f\in B_{X^*}$ so that
$\sum_{j=1}^q \alpha_j f(y^n_{i_j})\geq C_n$. Our construction gives
an $h\in B_{X^*}$ with $\sum_{j=1}^q \alpha_j h(y^n_{i_j})\geq C_n$
and $|h(y_j^{m})|<\epsilon_j$ if $m\leq q$ and $j\!\in\!\{m',\ldots
k\}\setminus\{i_1,\ldots i_q\}$. Because $(y_i^n)$ is bimonotone, we
may assume that $h(y_i^n)=0$ for all $i\geq n$ with $i> i_q$.  We
perturb $h$ by small multiples of the biorthogonal functionals of
$(y_i^n)$ to achieve $g\in X^*$ with $g(y_i^n)=h(y_i^n)$ for
$i\in\{i_1,\ldots,i_q\}$ and $g(y_i^n)=0$ for
$i\not\in\{i_1,\ldots,i_q\}$. Thus $g$ satisfies (\ref{E1}) and
(\ref{E2}).  All that remains is to check that $g\in
(1+\epsilon)B_{X^*}$.  Because $(y_i^n)$ is normalized and
bimonotone, we can estimate $\|g\|$ as follows:
$$\|g\|\leq\|h\|+\|g-h\|\leq1+\sum_{j=1}^{i_q-1}j\epsilon_j<1+\frac{\epsilon}{4}.
$$
\end{proof}

We are now prepared to give the proof of Proposition \ref{P2.4}. We
follow the same outline of the proof given in \cite{KO2} for
Proposition 3.4.

\begin {proof}[Proof of Proposition \ref{P2.4}]
Let $(x_i^n)$ be a normalized bad $uV$-array in $X$ and let $M_n$,
for $n\in\N$, be chosen so that the sequence $(x_i^n)_{i=n}^\infty$
is an $M_n$-bad $uV$-sequence and
$\lim_{n\rightarrow\infty}M_n=\infty$.  By Lemma \ref{L4.1} we just
need to consider the triangular array $(x_i^n)_{n\leq i}$. By
passing to a subarray using Lemma \ref{L4.2} and then renorming, we
may assume that $(x_i^n)_{n\leq i}$ is a normalized bimonotone basic
sequence in its lexicographical order.

We apply Theorem \ref{T4.3} for $\epsilon=1$ and $(C_n)=(M_n)$ to
obtain a subarray $(y_i^n)_{n\leq i}$ that satisfies the properties
(\ref{E1}) and (\ref{E2}). Moreover
 $(y^n_i)$ in its lexicographical order is a subsequence of
 $(x^n_i)$ in its lexicographical order, and thus is bimonotone.
 Furthermore, $(y^n_i)_{i=n}^\infty$ is a subsequence of
 $(x^n_i)_{i=n}^\infty$ for all $n\in\N$.  We denote
 $Y=[y_i^n]_{n\leq i}$.

Let $F(n)$ be a finite $\frac{1}{2n2^n}$-net in $[-2,2]$ which
contains the points 0,-2, and 2.  Whenever we have a functional
$g\in 2B_{X^*}$ which satisfies conditions (\ref{E1}) and (\ref{E2})
we may perturb $g$ by small multiples of the biorthogonal functions
of $(y_i^n)_{n\leq i}$ to obtain $f\in 3B_{X^*}$ which satisfies
(\ref{E1}), (\ref{E2}), and the following new condition
\begin{equation}\label{E3}
f(y_i^n)\in F(n)\qquad \textrm{for all } n\leq i.
\end{equation}
 We now start the construction of K.  Let $Y=[y^n_i]_{n\leq i}$ and
 $m\in N$.  We define the following,
$$L_m=\left\{(k_1,...,k_q)\;|\begin{array}{l}\; m\leq k_1< k_2<...<k_q,\\
||\sum_{i=1}^{q-1}\alpha_i y_{k_i}^m||\leq M_m \quad \textrm{for all
} (\alpha_i)\in B_V\\
||\sum_{i=1}^{q}\alpha_i y_{k_i}^m||> M_m \quad \textrm{for some }
(\alpha_i)\in B_V \end{array}\right\}
$$
It is important to note that if $(k_i)\in[\N]^\omega$ and $k_1\geq
m$ then there is a unique $q\in N$ such that $(k_1,...,k_q)\in L_m$.

Whenever $\vec{k}=(k_1,...,k_q)\in L_m$, our application of Theorem
\ref{T4.3} and then perturbation gives a functional $f\in 3B_{Y^*}$
which satisfies the properties (\ref{E1}),(\ref{E2}), and
(\ref{E3}). In particular we have that $\sum_{i=1}^q f(\alpha_i
y_{k_i}^m)>M_m$ for some $(\alpha_i)\in B_V$. We denote $f/3$ by
$f_{\vec{k}}$ and let for any $n\in\N$,
$$K_n=\{Q_m^* f_{\vec{k}}\;|\;m\in\N\;\vec{k}\in L_n\}.
$$
Here $Q_m$ denotes the natural projection of norm 1 from $Y$ onto
$[(y_i^n)]_{1\leq n\leq i\leq m}$. Finally, we define
$$K=\bigcup_{n=1}^\infty K_n\cup \{0\}.
$$
We first show that $(y_i^n|_K)_{n\leq i}$ is a bad u$V$-array as an
array in $C_b(K)$. Fix a column $n_0$. $(y_i^{n_0})_{i=n_0}^\infty$
is an $M_{n_0}$-bad u$V$-sequence. Consequently, given a subsequence
$(y_{k_i}^{n_0})_{i=1}^\infty$ of $(y_i^{n_0})_{i=n_0}^\infty$ we
have that $\vec{k}:=(k_1,...,k_q)\in L_{n_0}$ for some $q\in\N$. By
(\ref{E3}), $f_{\vec{k}}=Q^*_{q+1}f_{\vec{k}}$ and thus
$f_{\vec{k}}\in K_{n_0}\subset K$. $\sum_{i=1}^q
f_{\vec{k}}(\alpha_i y_{k_i}^{n_0})>\frac{M_{n_0}}{3}$ for some
$(\alpha_i)\in B_V$, and so we obtain that
$(y_i^{n_0}|_K)_{i=n_0}^\infty$ is an $(M_{n_0}/3)$-bad sequence in
$C_b(K)$, thus proving that $(y_i^n|_K)_{n\leq i}$ is a bad
u$V$-array.

$K$ is obviously a countable subset of $B_{Y^*}$. Since $Y$ is
separable, $K$ is $w^*$-metrizable.  Thus we need to show that $K$
is a $w^*$-closed subset of $B_{Y^*}$ in order to finish the proof.

Let $(g_j)\subset K$ and assume that $(g_j)$ converges $w^*$ to some
$g\in B_{Y^*}$.  We have to show that $g\in K$. Every $g_j$ is of
the form $Q^*_{m_j}f_{\vec{k}_j}$ for some $m_j\in\N$, $\vec{k}_j\in
L_{n_j}$, and some $n_j\in\N$.

By passing to a subsequence of $(g_j)$, we may assume that either
$n_j\rightarrow\infty$ as $j\rightarrow\infty$ or that there is an
$n\in N$ such that $n_j=n$ for all $j\in\N$.  We will start with the
first alternative. Let $i_j$ be the first element of $\vec{k}_j$.
Since $i_j\geq n_j$, we have that $i_j\rightarrow\infty$.  We also
have that $f_{\vec{k}_j}(y_i^n)=0$ for all $n\leq i< i_j$.  Thus
$f_{\vec{k}_j}\rightarrow0$ in the $w^*$ topology as
$j\rightarrow\infty$, so $g=0\in K$.

From now on we assume that there is an $n\in\N$ such that
$\vec{k}_j\in L_n$ for all $j\in\N$. $L_n$ is relatively
sequentially compact as a subspace of $\{0,1\}^\N$ endowed with the
product topology. Thus we may assume by passing to a subsequence of
$(g_j)$ that $\vec{k}_j\rightarrow \vec{k}$ for some $\vec{k}\in
\overline{L_n}$, the closure of $L_n$ in $\{0,1\}^\N$.

We now show that $\vec{k}$ is finite.  Suppose to the contrary that
$\vec{k}=(k_i)_{i=1}^\infty$. We have that $\vec{k}\in
\overline{L_n}$, so for all $r\in \N$ there exists $N_r\in\N$ such
that $\vec{k}_j=(k_1,...,k_r,\ell_1,...,\ell_s)$ for some
$\ell_1,...,\ell_s$ for all $j\geq N_r$.  Because $\vec{k}_j\in L_n$
we have that $k_1\geq n$, which implies that there exists $q\in\N$
such that $(k_1,...k_q)\in L_n$.  By uniqueness, $L_n$ does not
contain any sequence extending $(k_1,...,k_q)$. Therefore,
$\vec{k}_{N_{q+1}}=(k_1,...,k_{q+1}, \ell_1,...,\ell_s)\not\in L_n$,
a contradiction.

Since $B_{Y^*}$ is $w^*$-sequentially compact, we may assume that
$f_{\vec{k}_j}$ converges $w^*$ to some $f\in B_{Y^*}$.  We claim
that $f\in K$. To prove this we first show that $Q_m^*f\in K$ for
all $m\in\N$.  By (\ref{E2}) and (\ref{E3}) the set
$\{Q_m^*f_{\vec{k}_j}(y_i^n)\;|\; j\in\N\; 1\leq n\leq i\}$ has only
finitely many elements.  Since $Q_m^*f_{\vec{k}_j}\rightarrow
Q_m^*f$ as $j\rightarrow\infty$ we obtain that
$Q_m^*f_{\vec{k}_j}=Q_m^*f$ for $j\in\N$ large enough.  In
particular $Q_m^*f\in K$. Next let $q=\max \vec{k}$. Since
$\vec{k}_j\rightarrow\vec{k}$ and $\vec{k}$ is finite, we have
$Q_q^*f=f$ and thus $f\in K$.

Now we show that $g\in K$. By passing again to a subsequence of
$(g_j)$ we can assume that either $m_j\geq \max\vec{k}$ for all
$j\in\N$ or that there exists $m<\max\vec{k}$ such that $m_j=m$ for
all $j\in\N$.  If the first case occurs, then
$g_j=Q^*_{m_j}f_{\vec{k}_j}$ converges $w^*$ to $f$, and hence
$g=f\in K$. If the second case occurs then $g_j=Q^*_m f_{\vec{k}_j}$
converges $w^*$ to $Q_m^*f$, and hence $g=Q_m^*f\in K$.
\end{proof}

\section{Examples}
In previous sections, we introduced for any basic sequence $(v_i)$
the property $U_{(v_i)}$, and then proved that if a Banach space $X$
is $U_{(v_i)}$ then there exists a constant $C\geq1$ such that $X$
is $C-U_{(v_i)}$.  As Knaust and Odell proved that result for the
cases in which $(v_i)$ is the standard basis for $c_0$ or $\ell_p$
with $1\leq p<\infty$, we need to show that our result is not a
corollary of theirs.  For example, if $(v_i)$ is a basis for
$\ell_p\bigoplus\ell_q$ with $1<q<p<\infty$ which consists of the
union of the standard bases for $\ell_p$ and $\ell_q$ then  a Banach
space is $U_{(v_i)}$ or $C-U_{(v_i)}$ if and only if $X$ is
$U_{\ell_p}$ or $C-U_{\ell_p}$ respectively.  Thus the result for
this particular $(v_i)$ follows from \cite{KO2}. We make this idea
more formal by defining the following equivalence relation:

\begin{defn}
If $(v_i)$ and $(w_i)$ are normalized basic sequences then we write
$(v_i)\sim_U(w_i)$ (or $(v_i)\sim_{CU}(w_i)$) if each reflexive Banach space
is $U_{(v_i)}$ (or $C-U_{(v_i)}$) if and only if it is $U_{(w_i)}$
(or $C-U_{(w_i)}$).
\end{defn}

We define the equivalence relation strictly in terms of reflexive
spaces to avoid the unpleasant case of $\ell_1$. Because $\ell_1$
does  not contain any normalized weakly null sequence, $\ell_1$ is
trivially $U_{(v_i)}$ for every $(v_i)$. This is counter to the
spirit of what it means for a space to be $U_{(v_i)}$.  By
considering reflexive spaces, we avoid $\ell_1$, and we also make
the propositions included in this section formally stronger.
Reflexive spaces are also especially nice when considering
properties of weakly null sequences because the unit ball of a
reflexive spaces is weakly sequentially compact.  That is every
sequence in the unit ball of a reflexive space has a weakly
convergent subsequence.

In order to show that our result is not a corollary of the theorem
of Knaust and Odell, we give an example of a basic sequence $(v_i)$
such that $(v_i)\not\sim_U(e_i)$ where $(e_i)$ is the standard basis
for $c_0$ or $\ell_p$ with $1\leq p<\infty$. To this end we consider
a basis $(v_i)$ for a reflexive Banach space $X$ with the property
that $\ell_p$ is not $U_{(v_i)}$ for any $1<p<\infty$, but that $X$
is $U_{(v_i)}$ and not $U_{c_0}$.  We will be interested in
particular with the dual of the following space.
\begin{defn}
Tsirelson's space, T, is the completion of $\coo$ under the norm
satisfying the implicit relation:
$$||x||=||x||_\infty \vee\sup_{n\in\N (E_i)_1^n\subset[\N]^\omega
n\leq E_1 <...<E_n}\frac{1}{2}\sum_{i=1}^n ||E_i(x)||.
$$
$(t_i)$ is the unit vector basis of T and $(t^*_i)$ are the
biorthogonal functionals to $(t_i)$.
\end{defn}

Tsirelson constructed the dual of $T$ as the first example of a
Banach space which does not contain $c_0$ or $\ell_p$ for any $1\leq
p<\infty$ \cite{T}.  Though we are more interested in $T^*$, we use
the implicit definition of $T$ (which was formulated by Figiel and
Johnson in \cite{FJ}) as it is nice to work with.  Therefore, we
need some propositions that relate sequences in a space to sequences
in its dual.

\begin{prop}\label{P5.3}
If $(v_i)$ is a normalized bimonotone basic sequence, $(x_i)$ is a basic sequence, and $C>0$ then
\begin{enumerate}
\item[(i)] $(v_i)$ C-dominates $(x_i)$ if and only if $(v^*_i)$ is C-dominated by $(x^*_i)$,
\item[(ii)] $(v_i)$ C-dominates all of its normalized block bases if and only if $(v_i^*)$ is C-dominated by
all of its normalized block bases.
\end{enumerate}
\end{prop}
\begin{proof}
We assume that $(v_i)$ C-dominates $(x_i)$ and let $(a_i)\in\coo$ and $\epsilon>0$.
  There exists $(b_i)\in\coo$ such that $\sum a_i v_i^*(\sum b_i v_i)=||\sum a_i v_i^*||$
and $||\sum b_i v_i||<1+\epsilon$.  We have that
$$||\sum a_i v_i^*||=\sum a_i b_i=\sum a_i x_i^*(\sum b_i x_i)\leq C(1+\epsilon)||\sum a_i x_i^*||
.$$
Thus $(v_i^*)$ is $C-$dominated by $(x_i^*)$.  The converse is true by duality in the sense that we replace
 the roles of $(v_i)$ and $(x_i)$ by $(x_i^*)$
 and $(v_i^*)$ respectively. We have $(x_i^{**})=(x_i)$ and $(v_i^{**})=(v_i)$ and thus the converse follows and hence (i) is proven.

We assume that $(v_i)$ C-dominates all of its normalized block
bases, and let $(\sum_{j=k_i}^{k_{i+1}-1}a_j v^*_j)_{i=1}^\infty$ be
a normalized block basis of $(v_i)_{i=1}^\infty$. Because $(v_i)$ is
bimonotone, there exists a sequence $(b_j)_{j=1}^\infty\subset
[-1,1]$ such that: $$\sum_{j=k_i}^{k_{i+1}-1}a_j
v_j^*(\sum_{j=k_i}^{k_{i+1}-1}b_j
v_j)=1=||\sum_{j=k_i}^{k_{i+1}-1}b_j v_j|| \textrm{ for all }i\in\N.$$ Hence
$(\sum_{j=k_i}^{k_{i+1}-1}b_j v_j)_{i=1}^\infty$ is the sequence of
biorthogonal functions to the block basis $(\sum_{j=k_i}^{k_{i+1}-1}a_j
v^*_j)_{i=1}^\infty$. We have that $(v_i)$ C-dominates
$(\sum_{j=k_i}^{k_{i+1}-1}b_j v_j)$ and thus
$(\sum_{j=k_i}^{k_{i+1}-1}a_j v^*_j)$ C-dominates $(v_i^*)$ by (i).
The converse follows by duality and hence (ii) is proved.
\end{proof}

Proposition \ref{P5.3} together with some well known properties of
$(t_i)$ yields the following.

\begin{prop}\label{P5.4}
$(t^*_i)\not\sim_U(e_i)$ where $(e_i)$ is the standard basis for
$c_0$ or $\ell_p$ for $1\leq p<\infty$.
\end{prop}
\begin{proof}
It easily follows from the definition
that $(t_i)$ is an unconditional normalized basic sequence and that $(t_i)$ is dominated by each of
 its normalized block bases.  Also, the spreading model for $(t_i)$ is isomorphic to the standard $\ell_1$ basis. By
 proposition \ref{P5.3}, $(t_i^*)$ is an unconditional basic sequence that dominates all of its block bases and has its spreading model
 isomorphic to the standard basis for
 basis for $c_0$.  $T^*$ is reflexive because $(t_i^*)$ is unconditional and $T^*$ does not
  contain an isomorphic copy of $c_0$ or $\ell_1$.
   As $(t_i^*)$ has the standard basis for $c_0$ as its spreading model, we have that $\ell_p$ is
   not $U_{(t_i^*)}$ for all $1< p<\infty$.  Therefore $(t_i^*)\not\sim_U \ell_p$ for all $1\leq p<\infty$.
   As $(t_i^*)$ dominates all of its normalized block bases and every normalized weakly null sequence in $T^*$ has a subsequence
      equivalent to a normalized block basis of $(t_i^*)$, we have that $T^*$
   is $U_{(t_i^*)}$. $T^*$ does not contain $c_0$ isomorphically thus $T^*$ is not $U_{c_0}$.
   Therefore, $(t_i^*)\not\sim_U c_0$.
     \end{proof}

We have shown that $(t_i^*)\not\sim(e_i)$ where $(e_i)$ is the usual
basis for $c_0$ or $\ell_p$ for $1\leq p<\infty$,
 but we can actually say something much stronger than this.  One of the main properties of $\ell_p$ used in \cite{KO2} is
 that $\ell_p$ is subsymmetric.  If for each basic sequence $(v_i)$ there existed a constant $C\geq1$ and a subsymmetric
  basic sequence $(w_i)$ such that $(v_i)\sim_{CU}(w_i)$ then actually the first half of \cite{KO2} would apply to all
  basic sequences without changing anything. The following example shows in particular that this can not be done.

\begin{prop}\label{P5.5}
If $(v_i)$ is a normalized basic sequence such that $(v_{k_i})$
dominates $(v_i)$ for all $(k_i)\in[\N]^\omega$ then
$(v_i)\not\sim_U(t^*_i)$.
\end{prop}

In general, it can be fairly difficult to check if a Banach space is
$U_{(v_i)}$, as every normalized weakly null sequence in the space
needs to be checked.  In contrast to this, it is very easy to check
if $T^*$ is $U_{(v_i)}$.  This is because $(t_i)$ is dominated by
all of its block bases, and thus by Proposition \ref{P5.3} $T^*$ is
$U_{v_i}$ if and only if $(v_i)$ dominates a subsequence of
$(t^*_i)$.  In proving Proposition \ref{P5.5} we will carry this
idea further by considering a class of spaces, each of which have an
unconditional
 subsymmetric basis $(e_i)$ such that
$(e_i)$ is dominated by all of its normalized block bases.  The additional condition of subsymmetric gives that
$[e_i^*]$ is $U_{(v_i)}$ if and only if $(v_i)$ dominates $(e_i^*)$.  Hence, we need to check only one sequence instead
 of all weakly null sequences in $[e_i^*]$.

 The spaces we consider are generalizations of those
introduced by Schlumprecht \cite{S} as the first known arbitrarily
  distortable Banach spaces.
We put less restriction on the function $f$ given in the following
proposition, but we also infer less about the corresponding Banach
space.  The techniques used in \cite{S} are used to prove the
following proposition.

\begin{prop}\label{P5.6}
Let $f:\N\rightarrow[1,\infty)$ strictly increase to $\infty$,
$f(1)=1$, and $\lim_{n\rightarrow\infty}n/f(n)=\infty$. If $X$ is
defined as the closure of $\coo$ under the norm $||\cdot||$ which
satisfies the implicit relation:
$$||x||=||x||_{\infty}
\vee\sup_{m\in\N,E_1<...<E_m}\frac{1}{f(m)}\sum_{j=1}^m||E_j(x)||\quad\textrm{
for all }x\in\coo,
$$
then $X$ is reflexive.
\end{prop}

\begin{proof} Let $(e_n)$ denote the standard basis for $\coo$. It is straightforward to show that the norm
$||\cdot||$ as given in the statement of the theorem exists, as well
as that $(e_n)$ is a normalized, 1-subsymmetric and 1-unconditional
basis for $X$. Furthermore, $(e_n)$ is 1-dominated by all of its
normalized block bases.  We will prove that $X$ is reflexive by
showing that $(e_n)$ is boundedly complete and shrinking.

We first prove that $(e_n)$ is boundedly complete. As $(e_n)$ is
unconditional, if $(e_n)$ is not boundedly complete then it has some
normalized block basis which is equivalent to the standard $c_0$
basis.  However, $(e_n)$ is 1-dominated by all its normalized block
bases, so $(e_n)$ is also equivalent to the standard $c_0$ basis.
Hence $\sup_{N\in\N}||\sum_{n=1}^N e_n||<\infty$. This contradicts
that $||\sum_{n=1}^N e_n||\geq N/f(N)\rightarrow\infty$. Thus
$(e_n)$ is boundedly complete.

We now assume that $(e_n)$ is not shrinking. As $(e_n)$ is
unconditional, it has a normalized block basis $(x_n)$ which is
equivalent to the standard basis for $\ell_1$.  We will use James'
Blocking Lemma \cite{J} to show that this leads to a contradiction.
In one of its more basic forms, James' blocking lemma states that if
$(x_n)$ is equivalent to the standard basis for $\ell_1$ and
$\epsilon>0$ then $(x_n)$ has a normalized block basis which is
$(1+\epsilon)$-equivalent to the standard basis for $\ell_1.$
 Let
$0<\epsilon<\frac{1}{2}(f(2)-1)$.  By passing to a normalized block
basis using James' blocking lemma, we may assume that $(x_n)$ is
$(1+\epsilon)$-equivalent to the standard basis for $\ell_1$, and
thus any normalized block basis of $(x_n)$ will also be
$(1+\epsilon)$-equivalent to the standard basis for $\ell_1$. Let
$\epsilon_n>0$ such that $\sum_{n=1}^\infty\epsilon_n<\epsilon$.

We denote $||\cdot||_m$ to be the norm on $X$ which satisfies:
$$||x||_m=\sup_{E_1<...<E_m}\frac{1}{f(m)}\sum_{j=1}^m||E_j(x)||\quad\textrm{
for all }x\in\coo.
$$
We will construct by induction on $n\in\N$ a normalized block basis
$(y_i)$ of $(x_i)$ such that for all $m\in\N$ we have:
\begin{equation}\label{E4}
\textrm{If }||y_j||_m>\epsilon_j\textrm{ for some }1\leq j<n,
\textrm{ then }||y_n||_m<\frac{1+\epsilon_n}{f(m)}.
\end{equation}
For $n=1$ we let $y_1$=$x_1$, and note that (\ref{E4}) is vacuously
satisfied.

We now assume that we are given $n\geq1$ and finite block sequence
$(y_i)^n_{i=1}$ of $(x_i)$ which satisfies (\ref{E4}).  We have
$\lim_{m\rightarrow\infty}||y_i||_m\leq
\lim_{m\rightarrow\infty}\#\textrm{supp}(y_i)/f(m)=0$ (where
supp$(y_i)$ denotes the support of $y_i$). Thus, there exists $N>
\textrm{supp}(y_n)$ such that $||y_i||_m<\epsilon_i$ for all $1\leq
i\leq n$ and all $m\geq N$. Using James' blocking lemma, we block
$(x_i)_{i=N}^\infty$ into $(z_i)_{i=1}^\infty$ such that
$(z_i)_{i=1}^\infty$ is $(1+\epsilon_{n+1}/3)-$equivalent to the
standard $\ell_1$ basis. Let $M\geq6N/\epsilon_{n+1}$ and define
$y_{n+1}=\frac{1}{||\sum_{i=1}^{M}z_i||}\sum_{i=1}^{M}z_i$. Let
$m\in\N$ such that $||y_j||_m>\epsilon_j$ for some $1\leq j\leq n$.
By our choice of $N\in\N$, we have that $m<N$.  There exists
disjoint intervals $E_1,\ldots,E_m\subset\N$ and integers $1=k_0\leq
k_1\leq\ldots\leq k_m$ such that:
\begin{align*}
f(m)||y_{n+1}||_m&=\frac{1}{||\sum_1^M z_i||}\sum_{i=1}^m||P_{E_i}\sum_{j=k_{i-1}}^{k_i}z_j||\\
&\leq\frac{1+\epsilon_{n+1}/3}{M}\sum_{j=1}^m\left(||P_{E_i}z_{k_{i-1}}||+||\sum_{j=k_{i-1}+1}^{k_i-1}z_j||+||P_{E_i}z_{k_i}||\right)\\
&\leq\frac{1+\epsilon_{n+1}/3}{M}\left(M+2m\right)\\
&<\left(1+\epsilon_{n+1}/3\right)\left(1+2N/M\right)<\left(1+\epsilon_{n+1}/3\right)\left(1+\epsilon_{n+1}/3\right)<1+\epsilon_{n+1}.\\
\end{align*}
Hence, the induction hypothesis is satisfied.

We now show that property (\ref{E4}) leads to a contradiction with
$(y_i)$ being $(1+\epsilon)-$equivalent to the standard $\ell_1$
basis.  Let $n\in\N$.  We have for some $m\in\N$ that
$||\sum_{i=1}^n\frac{y_i}{n}||=||\sum_{i=1}^n\frac{y_i}{n}||_m$. By
(\ref{E4}) there exists $1\leq j\leq n+1$ such that
$||{y_i}||_m<\epsilon_i$ for all $1\leq i<j$ and
$f(m)||{y_i}||_m<1+\epsilon_i$ for all $j<i\leq n$. We have that:
\begin{align*}
||\sum_{i=1}^n\frac{y_i}{n}||=||\sum_{i=1}^n\frac{y_i}{n}||_m&\leq\frac{1}{n}\sum_{i=1}^{j-1}||y_i||_m+\frac{1}{n}||y_j||_m+\frac{1}{n}\sum_{i=j+1}^n
||y_i||_m\\
&< \frac{1}{n}\sum_{i=1}^{j-1}\epsilon_i +
\frac{1}{n}+\frac{1}{nf(m)}\sum_{i=j+1}^n
1+\epsilon_i\\
&<\frac{\epsilon}{n} + \frac{1}{n}+\frac{1}{f(2)}+\frac{\epsilon}{nf(2)}\\
&<\frac{\epsilon}{n} + \frac{1}{n}+\frac{1}{1+2\epsilon}+\frac{\epsilon}{n(1+2\epsilon))}.\\
\end{align*}
Thus $\inf_{n\in\N}||\sum_{i=1}^n
\frac{y_i}{n}||<\frac{1}{1+2\epsilon}$.  This contradicts that
$(y_i)$ is $(1+\epsilon)$ equivalent to the standard $\ell_1$ basis.
Hence $(e_i)$ is shrinking, and $X$ is reflexive.
\end{proof}

Using the reflexive spaces presented in Proposition \ref{P5.5}, we
can prove Proposition \ref{P5.5} for the case in which $(v_i)$ in
unconditional. The general case will then be reduced to this one.

\begin{lem}\label{L5.7}
If $(v_i)$ is a 1-suppression unconditional normalized basic sequence such that $(v_{k_i})$ dominates $(v_i)$ for all
$(k_i)\in[\N]^\omega$ and $(v_i)$ is not equivalent to the unit vector basis for
$c_0$ then there exists a reflexive Banach space $X$ which is $U_{(v_i)}$ and
not $U_{(t^*_i)}$.
\end{lem}
\begin{proof}
There exists $K\geq1$ such that $(v_{k_i})$ $K$-dominates $(v_i)$
for all $(k_i)\in[\N]^\omega$.  We define $\left<\cdot\right>$ to be
the norm on $(v_i)$ determined by:
$$\left<\sum_{i\in\N}a_i v^*_i\right>=\sup_{(k_i)\in[\N]^\omega}\left\Vert\sum_{i\in \N}a_i v^*_{k_i}\right\Vert \qquad\textrm{for all
}(a_i)\in\coo.
$$
Where $(v^*_i)$ is the sequence of biorthogonal functionals to
$(v_i)$.  The norm $\left<\cdot\right>$ is $K$-equivalent to the
original norm $||\cdot||$.  Furthermore, under the new norm
$(v_{k_i})$ $1$-dominates $(v_i)$ for all $(k_i)\in[\N]^\omega$.
Thus after possibly renorming, we may assume that $K$=1.

Let $\epsilon>0$ and $\epsilon_i\searrow0$ such that $\prod\frac{1}
{1-\epsilon_i}<1+\epsilon$. We have that $(v_i)$ is unconditional
and is not equivalent to the unit vector basis of $c_0$, so for all
$k\in\N$ there exists $N_{k}\geq k^2$ such that:
\begin{equation}
\left\Vert\sum
_{i=1}^{N_{k}}v_i\right\Vert>\frac{k+1}{\epsilon_{k+1}}.
\end{equation}

We define the function $f:\N\rightarrow [1,\infty)$ by:
$$f(n)= \begin{cases}
    1 &\text{if $n=1$},\\
    \frac{1}{1-\epsilon_1} &\text{ if $1<n\leq N_1$},\\
    k+1 &\text{if $N_{k}<n\leq N_{k+1}$ for $k\in\N$}.
          \end{cases}$$
We denote $|||\cdot|||$ to be the norm on $\coo$ determined by the
following implicit relation:

$$|||x|||=||x||_\infty\vee\sup_{m\in\N,E_1<...<E_m}\frac{1}{f(m)}\sum_{j=1}^m|||E_j(x)|||\quad\textrm{
for all }x\in\coo.
$$
The completion of $\coo$ under the norm $|||\cdot|||$ is denoted by
$X$, and its standard basis is denoted by $(e_i)$. We have that
$N_k>k^2$ which implies that $\lim_{k\rightarrow\infty}
k/f(k)=\infty$ and hence $X$ is reflexive by proposition \ref{P5.5}.

 We now show by induction on $k\in\N$ that if
$(a_i)_{i=1}^{N_{k}}\in\coo$ then
\begin{equation}\label{E5.1}
\left(\prod_{i=1}^k\frac{1}{1-\epsilon_i}\right)|||\sum_{i=1}^{N_k}
a_i e_i|||\geq||\sum_{i=1}^{N_k} a_i v^*_i||.
\end{equation}

For k=1, we have that $\frac{1}{1-\epsilon_1}|||\sum_{i=1}^{N_1}a_i
e_i|||\geq \sum_{i=1}^{N_1}|a_i|\geq||\sum_{i=1}^{N_1} a_i v^*_i||$.
Thus (\ref{E5.1}) is satisfied. Now we assume that $k\in\N$ and that
$(\prod_{i=1}^k\frac{1}{1-\epsilon_i})|||\sum_{i=1}^{N_k} a_i
e_i|||\geq||\sum_{i=1}^{N_k} a_i v^*_i||$ for all $(a_i)\in \coo$.

Let $(a_i)_{i=1}^{N_{k+1}}\subset\R$ such that
$||\sum_{i=1}^{N_{k+1}} a_i v^*_i||=1$. There exists
$(\beta_i)_{1=1}^{N_{k+1}}\subset\R$ such that $\sum\beta_i
a_i=||\sum\beta_i v_i||=1$.  Let $I=\left\{j\in\N:\quad
|\beta_j|<\frac{\epsilon_{k+1}}{k+1}\right\}$. If $\sum_{i\in
I}|a_i|\geq k+1$ then $|||\sum_{i=1}^{N_{k+1}} a_i
e_i|||\geq\frac{1}{k+1}\sum_{i\in I}|a_i|\geq 1=||\sum a_i v_i^*||$
and we are done.  Therefore we assume that $\sum_{i\in I}|a_i|<
k+1$, and thus $\sum_{i\in I}\beta_i a_i\leq\sum_{i\in
I}\frac{\epsilon_{k+1}}{k+1}|a_i|<\epsilon_{k+1}$. We let
$\{j_i\}_{i=1}^{\sharp I^C}=I^C$, and claim that $\sharp I^C\leq
N_{k}$. Indeed, if we assume to the contrary that $\sharp I^C>
N_{k}$, then

$$1\geq ||\sum_{i=1}^{\sharp
I^C}\beta_{j_i} v_{j_i}||\geq||\sum^{\sharp I^C}_{i=1}\beta_{j_i}
v_i||\geq\frac{\epsilon_{k+1}}{k+1}||\sum^{N_{k}}_{i=1}v_i||>\frac{\epsilon_{k+1}}{k+1}\frac{k+1}{\epsilon_{k+1}}=1
.$$ The first inequality is due to $(v_i)$ being 1-suppression
unconditional, and the second inequality is due to $(v_i)$ being
$1$-dominated by $(v_{j_i})$. Thus we have a contradiction and our
claim that $\sharp I^C\leq N_{k}$ is proven. We now have that
\begin{align*}
1=\sum\beta_i a_i&=\sum_I \beta_i a_i+\sum_{I^C}\beta_i a_i\\
&<\epsilon_{k+1}+||\sum_{i=1}^{\sharp I^C} a_{j_i} v^*_{j_i}||\\
&\leq\epsilon_{k+1}+||\sum_{i=1}^{\sharp I^C} a_{j_i} v^*_{i}||\\
&\leq\epsilon_{k+1}+(\prod_{i=1}^k\frac{1}{1-\epsilon_i})|||\sum_{i=1}^{\sharp I^C} a_{j_i} e_i|||\quad\textrm{by induction hypothesis}\\
&\leq\epsilon_{k+1}+(\prod_{i=1}^k\frac{1}{1-\epsilon_i})|||\sum_{i=1}^{N_{k+1}} a_i e_i|||\quad\textrm{by 1-subsymetric}.\\
\end{align*}
The last inequality gives that
  $1\leq(\prod_{i=1}^{k+1}\frac{1}{1-\epsilon_i})|||\sum_{i=1}^{N_{k+1}}
a_i e_i|||$.  Thus the induction hypothesis is satisfied.

 We have that $(e_i)$
dominates $(v^*_i)$, and hence $(v_i)$ dominates $(e^*_i)$.
$(e^*_i)$ is subsymmetric and dominates all its block bases, so
$[e_i^*]$ is $U_{(v_i)}$. $(e_i^*)$ is weakly null and is not
equivalent to the unit vector basis of $c_0$, so $[e_i^*]$ is not
$U_{(t_i^*)}$.
\end{proof}

The proof of Proposition \ref{P5.5} will follow easily once we have
established some properties of lower unconditional forms.

\begin{defn}
If $(v_i)$ is a normalized basic sequence, then the \textit{lower unconditional form} of $(v_i)$ is
the basic sequence $(u_i)$ determined by:
$$||\sum_{i\in\N}a_i u^*_i||=\sup_{E\subset\N}||\sum_{i\in E}a_i v^*_i|| \qquad\textrm{for all
}(a_i)\in\coo,
$$
where $(v_i^*)$ and $(u_i^*)$ are the sequences of biorthogonal functional to $(v_i)$ and $(u_i)$.
\end{defn}
The following lemma shows why the name lower unconditional form is appropriate and
also that some properties of basic sequences are passed on to their lower unconditional forms.

\begin{lem}\label{L5.9}
If $(v_i)$ and $(x_i)$ are normalized bimonotone basic sequences  with lower unconditional forms $(u_i)$ and $(y_i)$ respectively then
\begin{enumerate}
\item[(i)] $(u_i)$ is a 1-suppression unconditional normalized basic
sequence,
\item[(ii)] $(v_i)$ $1-$dominates $(u_i)$,
\item[(iii)] if $C\geq1$ and $(v_{k_i})$ $C-$dominates $(v_i)$ for all $(k_i)\in[\N]^\omega$ then $(u_{k_i})$ C-dominates $(u_i)$ for all
$(k_i)\in[\N]^\omega$,
\item[(iv)] if $(v_i)$ $C$-dominates $(x_i)$ then $(u_i)$ C-dominates $(y_i)$.
\end{enumerate}
\end{lem}
\begin{proof}
Let $(a_i)\in\coo$ and $J\subset\N$. We have that
$||u^*_i||=||v^*_i||=1$ because $(v_i)$ is normalized and
bimonotone. To check unconditionality we consider:
$$||\sum_{i\in\N}a_i u_i^*||=\sup_{E\subset\N}||\sum_{i\in E}a_i v_i^*||\geq\sup_{E\subset J}||\sum_{i\in E}a_i v_i^*||=||\sum_{i\in J}a_i
u_i^*||.
$$
Thus $(u_i^*)$ is a 1-suppression unconditional normalized basic sequence, and hence $(u_i)$ is as well.  Thus (i) is proven.

$(u_i^*)$ $1-$dominates $(v^*_i)$, and hence $(v_i)$ $1-$dominates
$(u_i)$ by Proposition \ref{P5.3} which proves (ii).

Let $(k_i)\in[\N]^\omega$ and assume that $(v_{k_i})$ $C$-dominates
$(v_i)$. $(v^*_i)$ $C-$dominates $(v^*_{k_i})$ by Proposition
\ref{P5.3}. Thus,
$$||\sum_{i\in\N}a_i u_i^*||=\sup_{E\subset\N}||\sum_{i\in E}a_i v_i^*||\leq C\sup_{E\subset\N}||\sum_{i\in E}a_i v_{k_i}^*||=C||\sum_{i\in\N}a_i
u_{k_i}^*||.
$$
Thus $(u^*_i)$ $C-$dominates $(u^*_{k_i})$, and hence (iii) follows
by Proposition \ref{P5.3}

Assume $(x_i^*)$ $C-$dominates $(v_i^*)$. Thus,
$$||\sum_{i\in\N}a_i u_i^*||=\sup_{E\subset\N}||\sum_{i\in E}a_i v_i^*||\leq C\sup_{E\subset\N}||\sum_{i\in E}a_i x_{i}^*||=||\sum_{i\in\N}a_i
y_{i}^*||.
$$
Hence, $(y_i^*)$ $C-$dominates $(u_i^*)$ and (iv) is proven.
\end{proof}
We now use lower unconditional forms to reduce the proof of
\ref{P5.5} to the case in Lemma \ref{L5.7}.
\begin{proof}[Proof of Proposition \ref{P5.5}]
Without loss of generality we assume that $(v_i)$ is bimonotone.
 We let $(u_i)$ be the lower unconditional form of $(v_i)$. If $T^*$ were not $U_{(v_i)}$ then $(v_i)\not\sim_U(t_i^*)$ because
$T^*$ is $U_{(t^*_i)}$.  Thus we may assume that $T^*$ is
$U_{(v_i)}$.  Hence, there exists $(k_i)\in[\N]^\omega$ such that
$(v_{i})$ dominates $(t^*_{k_i})$. We now have that $(u_i)$
dominates $(t^*_{k_i})$ by Lemma \ref{L5.9} since $(t^*_i)$ is
unconditional. Therefore $(u_i)$ is not isomorphic to the standard
basis for $c_0$.  We have by lemma \ref{L5.9} that $(u_{j_i})$
$1-$dominates $(u_i)$ for all $(j_i)\in[\N]^\omega$. By Lemma
\ref{L5.7} there exists a space $X$ which is $U_{(u_i)}$ but not
$U_{(t_i^*)}$. The space $X$ is also $U_{(v_i)}$ because $(v_i)$
dominates $(u_i)$, and thus our claim is proven.
\end{proof}
We also considered the question: "Does there exist a basic sequence
$(v_i)$ such that $(v_i)\not\sim_U(w_i)$ for any unconditional
$(w_i)$?". This is a much harder question, which is currently open.
Neither the summing basis for $c_0$, nor the standard basis for
James' space give a solution, as these are covered by the following
proposition:

\begin{prop}
If $(v_i)$ is a basic sequence such that
$\sup_{n\in\N}||\sum_{i=1}^n\epsilon_i v_i||<D$ for some
$(\epsilon_i)\in\{-1,1\}^\N$ and constant $D<\infty$ then
$(v_i)\sim_U c_0$.
\end{prop}
\begin{proof}
Let X be a C-$U_V$ Banach space, and let $(x_i)\in S_X$ be weakly
null. By Ramsey's theorem, we may assume by passing to a subsequence
that $(v_i)$ C-dominates every subsequence of $(x_i)$.  By a theorem
of John Elton \cite{E}, there exists $K<\infty$ and a subsequence
$(y_i)$ of $(x_i)$ such that if $(a_i)_1^\infty\in[-1,1]^\N$ and
$I\subset\{i:\quad |a_i|=1\}$ is finite then $||\sum_I a_i y_i||\leq
K\sup_{n\in\N}||\sum_{i=1}^n\epsilon_i y_i||$.  Thus we have for all
$A\in[\N]^{<\omega}$ that
$$||\sum_{i\in A} \epsilon_i y_i||\leq K\sup_{n\in\N}||\sum_{i=1}^n \epsilon_{i} y_i||\leq
KC\sup_{n\in\N}||\sum_{i=1}^n \epsilon_i v_i|| < KCD.$$

As this is true for all $A\in[N]^{<\omega}$, $(y_i)$ is equivalent
to the unit vector basis of $c_0$. Every normalized weakly null
sequence in X has a subsequence equivalent to $c_0$, so X is
$U_{c_0}$.
\end{proof}


\begin{thebibliography}{KO}
\bibitem[BP]{BP} C. Bessaga and A. Pelczy\'{n}ski, {\em Spaces of continuous functions
IV}, Studia Math., \textbf{19} (1960), 53-62.
\bibitem[E]{E} J.~Elton, {\em Weakly null normalized sequences in
Banach spaces}, Dissertation, Yale University, 1978.
\bibitem[FJ]{FJ} T. Figiel and W.B. Johnson, {\em A uniformly convex Banach space which contains
no $\ell_p$}, Compositio Math. \textbf{29} (1974), 179-190.
\bibitem[GM]{GM} W.T. Gowers and B. Maurey, {\em The unconditional basic sequence problem},
Annals of Mathematics \textbf{156} (2002), 797-833.
\bibitem[J]{J} R.C.~James,{\em Uniformly nonsquare Banach spaces},
Ann. of Math. (2)\textbf{80}(1964), 542-550.
\bibitem[JO]{JO} W.B.~Johnson and E.~Odell, {\em Subspaces of $L_p$ which embed into $\ell_p$}, Compos. Math.,
\textbf{28} (1974), 37-49.
\bibitem[KO1]{KO1} H.~Knaust and E.~Odell, {\em On $c_{0}$-sequences in Banach
spaces}, Isreal J. Math. \textbf{67} (1989), 153-169.
\bibitem[KO2]{KO2} H.~Knaust and E.~Odell, {\em Weakly null sequences with upper $\ell_p$
estimates}, pp.85-107, in: "Functional Analysis, Proceedings, The
University of Texas at Austin, 1987-89", E.~Odell, H.~Rosenthal
(eds.), Lecture Notes in Mathematics 1470, Springer-Verlag, Berlin
1991.
\bibitem[MR]{MR} B.~Maurey and H.~Rosenthal {\em Normalized weakly null sequences with no unconditional subsequences},
Studia Math., Volume 1200, 156 pp. Berlin and New Yotk: Springer-Verlag, (1977).
\bibitem[O]{O} E.~Odell, {\em Applications of Ramsey theorems to Banach space theory},
         Notes in Banach spaces, ed. H.E. Lacey, Univ. of Texas Press, Austin, TX
        (1980), 379--404.
\bibitem[S]{S} Th. Schlumprecht, {\em An arbitrarily distortable Banach space}, Israel J. Math.
    \textbf{76} (1991), 81-95.
\bibitem[T]{T} B.S.~Tsirelson, {\em Not every Banach space contains an embedding of $\ell_p$ or $c_0$},
Funct. Anal. Appl. \textbf{8} (1974), 138-141.
\end{thebibliography}
\end{document}